\documentclass{amsart}
\usepackage{graphicx}

\usepackage{amsmath}
\usepackage{amssymb}
\usepackage{amsthm}
\usepackage{calc}
\def\R{{\bf R}}

\def\K{{\ K}}
\def\Diff{\operatorname{Diff}}
\def\ln{\operatorname{ln}}

\usepackage[cp1251]{inputenc}
\usepackage[russian]{babel}

\begin{document}

\begin{center}
{
\bf About amenability of  subgroups of the group of
diffeomorphisms of the interval.\\}

 {
 \bf     E.T.Shavgulidze} \\
{\it Department of Mechanics and Mathematics, Moscow State
University, Moscow, 119991 Russia}
 \\
\end{center}

Averaging linear functional on the space continuous functions of
the group of diffeomorphisms of interval is found. Amenability of
several discrete subgroups of the group of diffeomorphisms
$\Diff^3([0,1])$ of interval is prove. In particular, a solution
of the problem of amenability of the Thompson's group $F$ is
given. See about the Thompson's group $F$ in [1].

  Let $\Diff^1_+([0,1])$  be the group of all diffeomorphisms of
class $C^1$ of interval $[0,1]$ that preserve the endpoints of
interval, and  let $\Diff^3_+ ([0,1])$  be the subgroup of
$\Diff^1_+([0,1])$ consisting of all diffeomorphisms  of class
$C^3([0,1])$, and
$$\Diff^3_0 ([0,1])=\{f \in \Diff^3_+ ([0,1]): f'(0)=f'(1)=1 \}.$$
The group
$\Diff^1_+([0,1])$ is equipped with the topology inherited from
the space $C^1([0,1])$.

For any positive  $\delta <1$, denote by $C_0^{1,\delta}([0,1])$
the set of all functions $f \in C^1([0,1])$ such that $ f(0)=0$
and $  \exists C>0 \ \ \forall t_1,t_2 \in [0,1] \ \
|f'(t_2)-f'(t_1)|<C |t_2-t_1|^\delta.$ Define a Banach structure
on the linear space $C_0^{1,\delta}([0,1])$  by a norm
$$\|f\|_{1,\delta}=|f'(0)|+
\sup \limits_{t_1,t_2 \in
[0,1]}\frac{|f'(t_2)-f'(t_1)|}{|t_2-t_1|^\delta}$$ for any
function $f \in C_0^{1,\delta}([0,1])$.

Let $\Diff^{1,\delta}_+([0,1])=\Diff^1_+([0,1])\bigcap
C_0^{1,\delta}([0,1]).$  It is easy to see that
$\Diff^{1,\delta}_+([0,1])$ is a subgroup of the group
$\Diff^1_+([0,1]).$ The subgroup $\Diff^{1,\delta}_+([0,1])$ is
equipped with the topology inherited from the space
$C_0^{1,\delta}([0,1])$.

Let $C_b(\Diff^{1,\delta}_+([0,1]))$ be the linear space of all
bounded continuous functions on the space
$\Diff^{1,\delta}_+([0,1])$, and let $C_b(\Diff^{1}_+([0,1]))$ be
the linear space of all bounded continuous functions on the space
$\Diff^{1}_+([0,1])$.

Introduce the functions $e_{1,\delta}:
\Diff^{1,\delta}_+([0,1])\to \R, \ \ e_{1,0}:\Diff^{1}_+([0,1])
\to \R$ by setting $e_{1,\delta}(g)=1$ for any $g\in
\Diff^{1,\delta}_+([0,1])$ and $e_{1,0}(f)=1$ \\ for any $f\in
\Diff^{1}_+([0,1])$. Let $F_g(f)=F(g^{-1}\circ f)$ for any $g \in
\Diff^3_0 ([0,1])$,\\ $f\in \Diff^{1,\delta}_+([0,1])$ and $F \in
C_b(\Diff^{1,\delta}_+([0,1]))$.

{\bf Theorem  1.\/} {\it For any positive $\delta <\frac{1}{2}$,
there exists a linear functional \\
$L_{\delta}: C_b(\Diff^{1,\delta}_+([0,1]))\to \R$ such that
$L_{\delta}(e_{1,\delta})=1$, $|L_{\delta} (F)|\leq \sup
\limits_{f \in \Diff^{1,\delta}_+([0,1])}|F(f)|$,  $L_{\delta}
(F)\geq 0$ for any nonnegative function $F \in
C_b(\Diff^{1,\delta}_+([0,1]))$, and \\
$L_{\delta}
(F_g)=L_{\delta} (F)$ for any $g \in  \Diff^3_0 ([0,1])$  and $F
\in C_b(\Diff^{1,\delta}_+([0,1]))$. \/}

The restriction of any function of the space $
C_b(\Diff^{1}_+([0,1]))$ on $\Diff^{1,\delta}_+([0,1])$ belon to
the space $ C_b(\Diff^{1,\delta}_+([0,1]))$. Hence we obtain the
following assertion.

{\bf Corollary  1.1.\/} {\it There exists a linear functional
$L_{0}: C_b(\Diff^{1}_+([0,1]))\to \R$ such that
$L_{0}(e_{1,0})=1$, $|L_{0} (F)|\leq \sup \limits_{f \in
\Diff^{1}_+([0,1])}|F(f)|$, $L_{0} (F)\geq 0$ for any nonnegative
function $F \in C_b(\Diff^{1}_+([0,1]))$, and $L_{0} (F_g)=L_{0}
(F)$ for any $g \in  \Diff^3_0 ([0,1])$ and $F \in
C_b(\Diff^{1}_+([0,1]))$. \/}

 We say that a discrete subgroup $G$ of
$\Diff^3_0([0,1])$ satisfies condition $(a)$, if  there is a such $C>0$ that \\
$\sup \limits_{t \in [0,1]}|\ln(g'_1(t))-\ln(g'_2(t))|\geq C$ for
any $g_1,g_2 \in G, g_1\neq g_2 $.

{\bf Theorem  2.\/} {\it If a discrete subgroup $G$ of
$\Diff^3_0([0,1])$ satisfies condition $(a)$, then the subgroup
$G$ is amenable. \/}

Proof.Let $B(G)$ be the linear space of all bounded  functions on
the group $G$.

Let positive $\delta <\frac{1}{2}$,  let
$$p(f)=|\ln(f'(0)|+
\sup \limits_{t_1,t_2 \in
[0,1]}\frac{|\ln(f'(t_2))-\ln(f'(t_1))|}{|t_2-t_1|^{\delta}}$$ and
$r(f)=\inf\limits_{h \in G}p(h^{-1}\circ f)$ for $f \in
\Diff^{1,\delta}_+([0,1])$, $\theta (t)=1-t$ for $t \in [0,1]$ and
$\theta (t)=0$ for $t>1$.

For any fixed $f \in \Diff^{1,\delta}_+([0,1])$,  $C>0$,   the set
of functions  $\{\psi : \psi (t)=\ln(g'(t)), \ g \in G, \ \ 
p_n(g\circ f)<C\} $ contain in a compact subset of the space
$C([0,1])$, therefore it is finite according to condition (а).
Hence, we can define the linear mapping $\pi_{\delta} : B(G) \to
C_b(\Diff^{1,\delta}_+([0,1]))$ by setting
$$\pi_\delta F(f)=\frac{\sum \limits_{h \in G}\theta (p_\delta (h^{-1}\circ f)-r(f))F(h)}
{\sum \limits_{h \in G}\theta (p_\delta (h^{-1}\circ f)-r(f))}.$$

Assign a linear functional $l: B(G)\to \R$ by setting
$l(F)=L_{\delta}(\pi_\delta F)$.

It is easy to see that
$$|l(F)|=|L_{\delta} (\pi_\delta F)|\leq \sup \limits_{f \in
\Diff^{1,\delta}_+([0,1])}|\pi_\delta F(f)|\leq \sup \limits_{g
\in G}|F(g)|,$$ $l(F)\geq 0$ for any nonnegative function $F \in
B(G)$, and $l(e_G)=1$, where $e_G(g)=1$ for all $g \in G$.

Denote by $F_g(h)=F(g^{-1} \circ h)$ for $F \in B(G)$,  $g,h \in
G$.

 We have
$$\pi_\delta F_g(f)=\frac{\sum \limits_{h \in G}\theta (p_\delta(h^{-1}\circ f)-r(f))F(g^{-1}\circ h)}
{\sum \limits_{h \in G}\theta (p_\delta(h^{-1}\circ f)-r(f))}=$$
$$=\frac{\sum \limits_{h \in G}\theta (p_\delta(h^{-1}\circ g\circ
f)-r(g\circ f))F(h)} {\sum \limits_{h \in G}\theta
(p_\delta(h^{-1}\circ g\circ f)-r(g\circ f))}= \pi_\delta F(g\circ
f),$$ hence $l(F_g)=L_{\delta}(\pi_\delta
F_g)=L_{\delta}(\pi_\delta F)=l(F)$,\\ which implies the assertion
of Theorem 2.

In [2] \`E.Ghys and V.Sergiescu proved that The Thompson's group
$F$ is isomorphic to a discrete subgroup $G$ of $\Diff^3_0([0,1])$
which satisfies condition $(a)$, and therefore we obtain the
following assertion.

{\bf Corollary  2.1.\/} {\it The Thompson's group $F$ is amenable.
\/}

For the proof Theorem  1 we need the following auxiliary
assertions.

{\bf Lemma  1.\/} {\it There exist positive constants  $c_1,c_2$
such that
$$c_1 2^{n+1}(n+1)!\leq\ \int\limits_{-\infty}\limits^{+\infty} ...\int\limits_{-\infty}\limits^{+\infty}
\frac{dx_1...dx_{n}}{\sqrt{(1+x_1^2)(1+(x_2 - x_1)^2) ...(1+(x_n -
x_{n-1})^2)(1+x_n^2)}}\leq$$
$$\leq c_2 2^{n+1}(n+1)!$$
 for any natural  $n$. \/}

Proof. The Fourier transform of the function
$\frac{1}{\sqrt{1+x^2}} $ is given by the Bessel\\ function of the
second kind
$$\frac{2}{\sqrt{2\pi}}\K_0 (y) =\frac{2}{\sqrt{2\pi}}\int\limits_{0}\limits^{+\infty}\frac{cos(xy)dx}{\sqrt{1+x^2}}
=\frac{1}{\sqrt{2\pi}}\int\limits_{-\infty}\limits^{+\infty}\frac{e^{-ixy}dx}{\sqrt{1+x^2}}$$
(see  [6]).

Hence we have
$$I= \int\limits_{-\infty}\limits^{+\infty} ...\int\limits_{-\infty}\limits^{+\infty}
\frac{dx_1...dx_{n}}{\sqrt{(1+x_1^2)(1+(x_2 - x_1)^2) ...(1+(x_n -
x_{n-1})^2)(1+x_n^2)}}=$$
$$=\frac{2^n}{\pi}\int\limits_{-\infty}\limits^{+\infty}(\K_0 (y))^{n+1}dx
=\frac{2^{n+1}}{\pi}\int\limits_{0}\limits^{+\infty}(\K_0
(y))^{n+1}dx.$$

There exist positive constants  $\varepsilon<1$,  $c$ such that
$-\ln\frac{y}{\varepsilon}<\K_0 (y)<-\ln\frac{y}{2}$ for any $y\in
(0,\varepsilon)$, and $0<\K_0 (y)<ce^{-y}$ for any $y\geq
\varepsilon$.

Thus,
$$I\leq \frac{2^{n+1}}{\pi}[ \int\limits_{0}\limits^{\varepsilon}
(-\ln(\frac{y}{2}))^{n+1}dy
+\int\limits_{\varepsilon}\limits^{+\infty} (ce^{-y})^{n+1}dy
]\leq$$
$$\leq \frac{2^{n+1}}{\pi}[ \int\limits_{0}\limits^{2}
(-\ln(\frac{y}{2}))^{n+1}dy +c^{n+1}
\int\limits_{\varepsilon}\limits^{+\infty} e^{-(n+1)y}dy ]=$$
$$=\frac{2^{n+1}}{\pi} [ 2 ((n+1)!)
+\frac{c^{n+1}e^{-(n+1)\varepsilon }  }{n+1} ],$$ and
$$I\geq  \frac{2^{n+1}}{\pi} \int\limits_{0}\limits^{\varepsilon}
(-\ln(\frac{y}{\varepsilon}))^{n+1}dy
=\frac{\varepsilon}{\pi}2^{n+1} (n+1)! , $$ which implies the
assertion of lemma 1.

Let
$$v_1(\tau)=\int\limits_{-\infty}\limits^{+\infty}\frac{d\tau_1}{\sqrt{(1+\tau_1^2)(1+(\tau-\tau_1)^2)}}$$
for $\tau \in \R$.

{\bf Lemma  2.\/} {\it  The derivation $v_1'(t)$ is negative for
any $t>0$, and  $|v_1'(t)|\leq \frac{4}{|t|}v_1(t)$ for any $t\neq
0$. \/}

Proof. The derivation of the function $v_1(t)$ is equal
$$v_1'(t)=-\int\limits_{-\infty}\limits^{+\infty}\frac{(t-\tau)d\tau}{\sqrt{(1+\tau^2)(1+(t-\tau)^2)^3}}.$$
Taking $\tau_1= t-\tau$ we receive
$$v_1'(t)=-\int\limits_{-\infty}\limits^{+\infty}\frac{\tau_1 d\tau_1}{\sqrt{(1+\tau_1^2)^3(1+(\tau_1-t)^2)}}=$$
$$=-\int\limits_{0}\limits^{+\infty}(\frac{1}{\sqrt{1+(\tau_1-t)^2}}-\frac{1}{\sqrt{1+(\tau_1+t)^2}})\frac{\tau_1 d\tau_1}{\sqrt{(1+\tau_1^2)^3}}=$$
$$=-\int\limits_{0}\limits^{+\infty} \frac{4 \ t \ \tau_1 ^2 \
d\tau_1}{(\sqrt{1+(\tau_1-t)^2}+\sqrt{1+(\tau_1+t)^2})\sqrt{(1+\tau_1^2)^3(1+(\tau_1-t)^2)(1+(\tau_1+t)^2)}}.$$
Hence, if $t>0$, then $v_1'(t)<0$.

If $t\neq 0$, then
$$|v_1'(t)| \leq \int\limits_{0}\limits^{+\infty} \frac{4 \ |t| \ \tau_1 ^2 \
d\tau_1}{(1+\tau_1^2)(1+(\tau_1+|t|)^2)\sqrt{(1+(\tau_1-t)^2)(1+\tau_1^2)}}
\leq $$
$$  \leq \frac {4} { |t|} \int\limits_{0}\limits^{+\infty} \frac{
d\tau_1}{\sqrt{(1+(\tau_1-t)^2)(1+\tau_1^2)}}=\frac
{4}{|t|}v_1(t),$$ which implies the assertion of lemma 2.

 Let $v(t)=v_1(\ln(t+\sqrt{t ^2-1})$ for $ t\geq 1$.

 Let us denote
$ D_n=\{(x_1,...,x_{n-1}):0<x_1<...<x_{n-1}<1\}$ \\ and
$x_{-1}=x_{n-1}-1, \  x_0=0, \ x_n=1$.

Let
$$u_{1,n}(x_1,...,x_{n-1})=
\prod \limits_{k=1}
\limits^{n}\frac{1}{x_k-x_{k-1}}v(\frac{x_k-x_{k-2}}{2\sqrt{(x_k-x_{k-1})(x_{k-1}-x_{k-2})}}),$$
$$J_n=\int\limits_{0}\limits^{1} \int\limits_{x_1}\limits^{1}...\int\limits_{x_{n-2}}\limits^{1}
u_{1,n}(x_1,x_2,...,x_{n-1})dx_1dx_2...dx_{n-1},$$
$u_n(x_1,...,x_{n-1})=\frac{1}{J_n}u_{1,n}(x_1,...,x_{n-1})$.

{\bf Lemma  3.\/} {\it  The following equality is fulfilled
$$J_n=2^{n-1}\int\limits_{-\infty}\limits^{+\infty} ...\int\limits_{-\infty}\limits^{+\infty}
\frac{dt_1 dt_2...dt_{2n-1}}{\sqrt{(1+t_1^2)(1+(t_2 - t_1)^2)
...(1+(t_{2n-1} - t_{2n-2})^2)(1+t_{2n-1}^2)}}
$$
for any natural $n $, and $c_1 2^{3n-1}(2n)! \leq J_{n}\leq c_2
2^{3n-1}(2n)!$ . \/}

Proof. Taking $y_k=\frac{x_k-x_{k-1}}{1-x_{n-1}}$ we have
$$J_n=\int\limits_{0}\limits^{1} \int\limits_{x_1}\limits^{1}...\int\limits_{x_{n-2}}\limits^{1}
u_{1,n}(x_1,x_2,...,x_{n-1})dx_1dx_2...dx_{n-1}=$$
$$=\int\limits_{0}\limits^{+\infty} \int\limits_{0}\limits^{+\infty}...\int\limits_{0}\limits^{+\infty}
v(\frac{y_{1}+1}{2\sqrt{y_{1}}})
v(\frac{1+y_{n-1}}{2\sqrt{y_{n-1}}})\prod \limits_{k=2}
\limits^{n-1}v(\frac{y_k+y_{k-1}}{2\sqrt{y_k y_{k-1}}})
\frac{dy_1dy_2...dy_{n-1}}{y_1y_2...y_{n-1}}.$$

Making substitutions $t_{2k}=\frac{1}{2}\ln(y_k)$ and taking into
account
$$\ln( \frac{y_k+y_{k-1}}{2\sqrt{y_k y_{k-1}}} +
\sqrt{(\frac{y_k+y_{k-1}}{2\sqrt{y_k y_{k-1}}})^{2} -1}\, ) =
|t_{2k}-t_{2k-2}|,$$
$$v(\frac{y_k+y_{k-1}}{2\sqrt{y_k
y_{k-1}}})=v_1(|t_{2k}-t_{2k-2}|)=$$
$$=2\int\limits_{-\infty}\limits^{+\infty} \frac{d t_{2k-1}}{\sqrt{(1+ (t_{2k-1}-t_{2k-2}^2)(1+ (t_{2k}-t_{2k-1})^2)}}, $$
$$v(\frac{y_1+1}{2\sqrt{y_1}})=2\int\limits_{-\infty}\limits^{+\infty} \frac{d t_{1}}{\sqrt{(1+ t_{1}^2)(1+ (t_{2}-t_{1})^2)}}, $$
$$v(\frac{1+y_{n-1}}{2\sqrt{y_{n-1}}})=2\int\limits_{-\infty}\limits^{+\infty} \frac{d t_{2n-1}}{\sqrt{(1+ t_{2n-1}^2)(1+ (t_{2n-1}-t_{2n-2})^2)}}, $$
we receive
$$J_n=2^{n-1}\int\limits_{-\infty}\limits^{+\infty} ...\int\limits_{-\infty}\limits^{+\infty}
\frac{dt_1 dt_2...dt_{2n-1}}{\sqrt{(1+t_1^2)(1+(t_2 - t_1)^2)
...(1+(t_{2n-1} - t_{2n-2})^2)(1+t_{2n-1}^2)}}.$$

  Thus   the  assertion  of  lemma 3 is an immediate  consequence  of
lemma 1.

{\bf Lemma  4.\/} {\it  For arbitrary positive $\varepsilon $, it
exist positive constants $c_3$  such that
$$v(\frac{y_{1}+a}{2\sqrt{y_{1}a}})
v(\frac{a+y_{2}}{2\sqrt{ay_{2}}}) \leq c_3
v(\frac{y_{1}+y_{2}}{2\sqrt{y_{1}y_{2}}})$$ for all $a,y_1,y_2$,
satisfying  $\varepsilon \leq a \leq 1 ,y_1>0 ,y_2>0, y_1 +
y_2\leq 1$. \/}

Proof. Let
$t_{1}=-\frac{1}{2}\ln(y_1),t_{2}=-\frac{1}{2}\ln(y_2),\alpha=-\frac{1}{2}\ln(a)$
then $0\leq \tau_0 \leq r=-\frac{1}{2}\ln(\varepsilon)$ and
$$v(\frac{y_{1}+a}{2\sqrt{y_{1}a}})=v_1(t_1-\alpha ),
v(\frac{a+y_{2}}{2\sqrt{ay_{2}}}) =v_1(t_2-\alpha ),
v(\frac{y_{1}+y_{2}}{2\sqrt{y_{1}y_{2}}}) v_1(t_2-t_1 ).$$

It follows from the lemma 2 that the derivation $v_1'(t)$ is
negative for any $t>0$ and the function $v_1(t)$ decrease on
$[0,+\infty)$.

Thus, $v_1(t)\leq v_1(0)
=\int\limits_{-\infty}\limits^{+\infty}\frac{
d\tau}{(1+\tau^2)}=\pi$ for any $t\geq 0$.

As $v_1(t)>0$ for any $t \in \R$, the function $v_1(t)$ is
continuous and
$$\lim \limits_{t \to +\infty} \frac{v_1(t-r)}{v_1(t)}=1,$$ then it exists
positive constant $c^*$  such that $v_1(t-\tau)\leq c^* v_1(t)$
for any \\ $t>0, \tau \in[0,r]$.

Let  $ c_3=\pi(c^*)^2 .$

If $t_2\geq t_1$, then $t_2\geq t_2-t_1 \geq 0$  and
$$v_1(t_1-\alpha)v_1(t_2-\alpha) \leq (c^*)^2 v_1(t_1)v_1(t_2)
\leq \pi (c^*)^2 v_1(t_2) \leq \pi(c^*)^2  v_1(t_2-t_1)= c_3
v_1(t_2-t_1).$$

If $t_1\geq t_2$, then $v_1(t_1-\alpha)v_1(t_2-\alpha) \leq c_3
v_1(t_1-t_2)= c_3 v_1(t_2-t_1).$

Hence, $v(\frac{y_{1}+a}{2\sqrt{y_{1}a}})
v(\frac{a+y_{2}}{2\sqrt{ay_{2}}}) \leq c_3
v(\frac{y_{1}+y_{2}}{2\sqrt{y_{1}y_{2}}})$.

Let us define a function $\vartheta (t)=1$  for $t \in [0,1]$ and
$\vartheta (t)=0$ for $t \in (-\infty ,0)\bigcup (1,+\infty)$.

{\bf Lemma 5.\/} {\it  The following equality is valid
$$\lim \limits_{n \to \infty} (\int\limits_{0}\limits^{1} \int\limits_{x_1}\limits^{1}...\int\limits_{x_{n-2}}\limits^{1}
(1-\vartheta (\frac{1}{\varepsilon}\max \limits_{1\leq k \leq n}
(x_k-x_{k-1}))) u_{n}(x_1,x_2,...,x_{n-1})dx_1dx_2...dx_{n-1})=0$$
for any positive $\varepsilon <1$. \/}

Proof. We have
$$\int\limits_{0}\limits^{1} \int\limits_{x_1}\limits^{1}...\int\limits_{x_{n-2}}\limits^{1}
(1-\vartheta (\frac{1}{\varepsilon}\max \limits_{1\leq k \leq n}
(x_k-x_{k-1})))
u_{n}(x_1,x_2,...,x_{n-1})dx_1dx_2...dx_{n-1}\leq$$
$$\leq \sum \limits_{k=1}\limits^{n}
\int\limits_{0}\limits^{1}
\int\limits_{x_1}\limits^{1}...\int\limits_{x_{n-2}}\limits^{1}
(1-\vartheta (\frac{1}{\varepsilon} (x_k-x_{k-1})))
u_{n}(x_1,x_2,...,x_{n-1})dx_1dx_2...dx_{n-1}.$$

Let $1\leq k \leq n$. Let us  take $r=x_{k}-x_{k-1} \ \
(\varepsilon\leq r \leq 1)$,
$$y'_{-1}=x_{n-1}-1, \ y'_0=0, \ y'_1=x_{1},..., \ y'_{k-1}=x_{k-1},$$
$$y'_{k}=x_{k+1}-r,...,\
  y'_{n-2}=x_{n-1}-r, \ y'_{n-1}=1-r$$
and $y_l=\frac{y'_l}{1-r}$ for $l=-1,0,1,...,n-1$.

It follows from lemma 4 that it exists positive $c_3$ such that
$$v(\frac{x_{k}-x_{k-2}}{2\sqrt{(x_{k}-x_{k-1})(x_{k-1}-x_{k-2})}})
v(\frac{x_{k+1}-x_{k-1}}{2\sqrt{(x_{k+1}-x_{k})(x_{k}-x_{k-1})}})\leq$$
$$\leq c_3
v(\frac{x_{k+1}-x_{k}+x_{k-1}-x_{k-2}}{2\sqrt{(x_{k+1}-x_{k})(x_{k-1}-x_{k-2})}})=c_3
v(\frac{y'_{k}-y'_{k-2}}{2\sqrt{(y'_{k}-y'_{k-1})(y'_{k-1}-y'_{k-2})}})$$.

Hence,
$$I_k=\int\limits_{0}\limits^{1}
\int\limits_{x_1}\limits^{1}...\int\limits_{x_{n-2}}\limits^{1}
(1-\vartheta (\frac{1}{\varepsilon} (x_k-x_{k-1})))
u_{n}(x_1,x_2,...,x_{n-1})dx_1dx_2...dx_{n-1}\leq$$
$$\leq \frac{c_3}{J_n}\int\limits_{\varepsilon}\limits^{1} \frac{1}{r} [\int\limits_{0}\limits^{1-r}
\int\limits_{y'_1}\limits^{1-r}...\int\limits_{y'_{n-3}}\limits^{1-r}
u_{1,n-1}(y'_1,y'_2,...,y'_{n-1})dy'_1dy'_2...dy'_{n-2}]dr=$$
$$=\frac{c_3}{J_n}\int\limits_{\varepsilon}\limits^{1} \frac{(1-r)dr}{r} \cdot \int\limits_{0}\limits^{1}
\int\limits_{y_1}\limits^{1}...\int\limits_{y_{n-3}}\limits^{1}
u_{1,n-1}(y_1,y_2,...,y_{n-1})dy_1dy_2...dy_{n-2}\leq $$
$$\leq |\ln \varepsilon | \frac{c_3 J_{n-1}}{J_n} \leq |\ln \varepsilon |
\frac{c_3 c_{2}}{16 c_{1}(2n-1)n} .$$

Thus,
$$0\leq \int\limits_{0}\limits^{1} \int\limits_{x_1}\limits^{1}...\int\limits_{x_{n-2}}\limits^{1}
(1-\vartheta (\frac{1}{\varepsilon}\max \limits_{1\leq k \leq n}
(x_k-x_{k-1})))
u_{n}(x_1,x_2,...,x_{n-1})dx_1dx_2...dx_{n-1}\leq$$
$$ \leq n|\ln \varepsilon |
\frac{c_3 c_{2}}{16 c_{1}(2n-1)n}=|\ln \varepsilon | \frac{c_3
c_{2}}{16 c_{1}(2n-1)} ,$$ which implies the assertion of lemma 5.

{\bf Lemma  6.\/} {\it  The following equality is fulfilled
$$\lim \limits_{n \to \infty} (\int\limits_{0}\limits^{1} \int\limits_{x_1}\limits^{1}...\int\limits_{x_{n-2}}\limits^{1}
\vartheta (\frac{1}{r}\min \limits_{1\leq k \leq n}
\frac{x_{k}-x_{k-2}}{2\sqrt{(x_{k}-x_{k-1})(x_{k-1}-x_{k-2})}})$$
$$u_{n}(x_1,x_2,...,x_{n-1})dx_1dx_2...dx_{n-1})=0$$ for any  $r>1$. \/}

Proof. We have
$$ \int\limits_{0}\limits^{1} \int\limits_{x_1}\limits^{1}...\int\limits_{x_{n-2}}\limits^{1}
\vartheta (\frac{1}{r}\min \limits_{1\leq k \leq n}
\frac{x_{k}-x_{k-2}}{2\sqrt{(x_{k}-x_{k-1})(x_{k-1}-x_{k-2})}})$$
$$u_{n}(x_1,x_2,...,x_{n-1})dx_1dx_2...dx_{n-1}\leq$$
$$ \leq \sum \limits_{k=1}\limits^{n}
\int\limits_{0}\limits^{1}
\int\limits_{x_1}\limits^{1}...\int\limits_{x_{n-2}}\limits^{1}
\vartheta (\frac{1}{r}
\frac{x_{k}-x_{k-2}}{2\sqrt{(x_{k}-x_{k-1})(x_{k-1}-x_{k-2})}})$$
$$u_{n}(x_1,x_2,...,x_{n-1})dx_1dx_2...dx_{n-1}.$$

Let $1\leq k \leq n$. Making the same substitutions
$y_l=\frac{x_l-x_{l-1}}{1-x_{n-1}}$ и $t_{2l}=\frac{1}{2}\ln(y_l)$
as at lemma 3 we receive
$$ I_k=
\int\limits_{0}\limits^{1}
\int\limits_{x_1}\limits^{1}...\int\limits_{x_{n-2}}\limits^{1}
\vartheta (\frac{1}{r}
\frac{x_{k}-x_{k-2}}{2\sqrt{(x_{k}-x_{k-1})(x_{k-1}-x_{k-2})}})$$
$$u_{n}(x_1,x_2,...,x_{n-1})dx_1dx_2...dx_{n-1}=$$
$$=\frac{1}{J_n}\int\limits_{-\infty}\limits^{+\infty} ...\int\limits_{-\infty}\limits^{+\infty}
\vartheta (\frac{1}{r}
\frac{e^{2t_{2k}}+e^{2t_{2(k-1)}}}{2e^{t_{2k}}e^{t_{2(k-1)}}})$$
$$\frac{dt_1 dt_2...dt_{2n-1}}{\sqrt{(1+t_1^2)(1+(t_2 - t_1)^2)
...(1+(t_{2n-1} - t_{2n-2})^2)(1+t_{2n-1}^2)}}.$$

The inequality
$\frac{e^{2t_{2k}}+e^{2t_{2(k-1)}}}{2e^{t_{2k}}e^{t_{2(k-1)}}}\leq
r$ is valid if and only if
$$|t_{2k} - t_{2k-2}|\leq a=\ln(r+\sqrt{r ^2-1}).$$

Hence, $ \vartheta (\frac{1}{r}
\frac{e^{2t_{2k}}+e^{2t_{2(k-1)}}}{2e^{t_{2k}}e^{t_{2(k-1)}}})=
\vartheta (\frac{1}{a}|t_{2k} - t_{2k-2}|) $ и
$$I_k =\frac{2^{n-1}}{J_n}\int\limits_{-\infty}\limits^{+\infty} ...\int\limits_{-\infty}\limits^{+\infty}
\frac{\vartheta (\frac{1}{a}|t_{2k} - t_{2k-2}|)dt_1
dt_2...dt_{2n-1}}{\sqrt{(1+t_1^2)(1+(t_2 - t_1)^2) ...(1+(t_{2n-1}
- t_{2n-2})^2)(1+t_{2n-1}^2)}}.$$

As $|t_{2k} - t_{2k-2}|\leq a$ we have $|t_{2k}-t_{2k-3}|\leq
|t_{2k-2}-t_{2k-3}|+a$. Thus
$$1+(t_{2k}-t_{2k-3})^2\leq (1+(t_{2k-2}-t_{2k-3})^2)4(1+a)^2$$
and
$$ \frac{1}{1+(t_{2k-2}-t_{2k-3})^2}\leq
\frac{4(1+a)^2}{1+(t_{2k}-t_{2k-3})^2}.$$

Taking into account
$$\int\limits_{-\infty}\limits^{+\infty}
 \frac{dt_{2k-1}}{\sqrt{(1+(t_{2k} -
t_{2k-1})^2)(1+(t_{2k-1}-t_{2k-2})^2)}}=v_1(|t_{2k}-t_{2k-2}|)\leq
\pi$$ we receive
$$I_k \leq \frac{\pi 2^{n}(1+a)}{J_n}\int\limits_{-\infty}\limits^{+\infty} ...\int\limits_{-\infty}\limits^{+\infty}
\int\limits_{-\infty}\limits^{+\infty}
...\int\limits_{-\infty}\limits^{+\infty}
(\int\limits_{t_{2k}-a}\limits^{t_{2k}+a}dt_{2k-2})$$
$$
\frac{dt_1
dt_2...dt_{2k-3}dt_{2k}...dt_{2n-1}}{\sqrt{(1+t_1^2)(1+(t_2 -
t_1)^2) ...(1+(t_{2k-3} - t_{2k})^2)...(1+(t_{2n-1} -
t_{2n-2})^2)(1+t_{2n-1}^2)}}=$$
$$= \frac{8\pi a(1+a)J_{n-1}}{J_n} \leq \frac{c_2 \pi a(1+a)}{
c_{1}(2n-1)n} .$$

Hence
$$0\leq \int\limits_{0}\limits^{1} \int\limits_{x_1}\limits^{1}...\int\limits_{x_{n-2}}\limits^{1}
\vartheta (\frac{1}{r}\min \limits_{1\leq k \leq n}
\frac{x_{k}-x_{k-2}}{2\sqrt{(x_{k}-x_{k-1})(x_{k-1}-x_{k-2})}})$$
$$u_{n}(x_1,x_2,...,x_{n-1})dx_1dx_2...dx_{n-1}\leq$$
$$ \leq n
\frac{c_2 \pi a(1+a)}{ c_{1}(2n-1)n}= \frac{c_2 \pi a(1+a)}{
c_{1}(2n-1)} ,$$ which implies the assertion of Lemma 6.

Define the mapping
 $A:\Diff^1_+([0,1])\to C_0([0,1])$ by setting
$$
A(q)(t)=\ln(f'(t))-\ln(q'(0)) \qquad \forall t\in[0,1].
$$

The mapping $A$ is a topological isomorphism between the space
$\Diff^1_+([0,1])$, $C_0([0,1])$ moreover
$$
A^{-1} (\xi) (t)=\frac {\int_0^t e^{\xi (\tau)}d\tau} {\int_0^1
e^{\xi (\tau)}d\tau}.
$$

 Introduce the Wiener measure $w$ on the space $C_0([0,1])$.
Define a Borel  measure $\nu$  on $\Diff^1_+([0,1])$ by setting
$\nu(X)=w(A(X))$ for any Borel subset $X$ of topological space
$\Diff^1_+([0,1])$.

Let $ \delta \in (0, \frac{1}{2})$. It follows from the properties
of Wiener measure $w$ (see [5]) that measure  $\nu$ is
concentrated on the set $E_\delta=\Diff^{1,\delta}_+([0,1])$, i.d.
$\nu(E_\delta)=1$, moreover the Borel subsets of metric space
$E_\delta$ is measurable with respect to the measure $\nu$.

As was proved in [3], the  measure $\nu$ is quasi-invariant with
respect to the left action of subgroup $\Diff^3_+([0,1])$ on the
group $\Diff^1_+([0,1])$, moreover
\begin{eqnarray}
\nu (g  X)\,=\, \frac 1 {\sqrt{g '(0) g '(1)}} \int_X \,e^{\frac
{g ''(0)} {g '(0)} f'(0)- \frac {g ''(1)} {g '(1)} f'(1) +\int_0^1
S_g (q(t)) (q'(t))^2 dt} \,\nu (d q), \nonumber
\end{eqnarray}
for any  Borel subset $X$  of topological space $\Diff^1_+
([0,1])$,and any
$g \in \Diff^3_+ ([0,1])$,  \\
where $gX=\{ g\circ q : q \in X \}$  and  $ S_g (\tau) = \frac {g
'''(\tau)} {g '(\tau)}- \frac 3 2 (\frac {g ''(\tau)} {g '(\tau)
})^2$ (Schwartz derivative of function  $g$).

{\bf Lemma  7.\/} {\it The following equality is valid
$$\int\limits_{E_\delta}(q'(0))^l\,\nu (d q)=
 \int\limits_{E_\delta}(q'(1))^l\,\nu (d q)$$ for any natural $l$.\/}

Proof. Let $\xi = A(q)$, i.d. $\xi(t)=\ln(f'(t))-\ln(q'(0)) $.
Then
$$
q' (0)=\frac {1} {\int_0^1 e^{\xi (\tau)}d\tau}, \  q' (1)=\frac
{e^{\xi (1)}} {\int_0^1 e^{\xi (\tau)}d\tau}.
$$

Let us take
$$
M_l= \int\limits_{E_\delta}(q'(1))^l\,\nu (d q)
=\int\limits_{\Diff^1_+([0,1])}(q'(1))^l\,\nu (dq)=$$
$$=\int\limits_{C_0([0,1])}(\frac {e^{\xi (1)}} {\int_0^1
e^{\xi (\tau)}d\tau})^l\,\nu (d q)=
 \int\limits_{C_0([0,1])}(\frac {1} {\int_0^1 e^{\xi
(1-\tau)-\xi (1)}d\tau})^l\,\nu (d q)$$

Let $\zeta (t)=\xi (1-t)-\xi (1)$. The Wiener measure $w$ is
invariant with respect to the action $\zeta \longmapsto \xi $,
thus,
$$
M_l= \int\limits_{C_0([0,1])}(\frac {1} {\int_0^1 e^{\zeta
(\tau)}d\tau})^l\,\nu (d q)=$$
$$=\int\limits_{\Diff^1_+([0,1])}(q'(0))^l\,\nu (dq)=
\int\limits_{E_\delta}(q'(0))^l\,\nu (d q),$$ which implies the
assertion of Lemma 7.

  Introduce the measure $\nu _n=\nu \otimes ...\otimes\nu $ on the space
$E_\delta^n=E_\delta \times ...\times E_\delta $.

Let $c_4=1+M_1+M_2+\int\limits_{E_1}\int_0^1  (q'(t))^2 dt \,\nu
(d q)$.

For any $r>0$,  $g \in \Diff^3_+([0,1])$,
$\overline{x}=(x_1,...,x_{n-1}) \in D_n$, we write \\
$C_g=1+\max \limits_{0\leq t \leq 1}(|\frac {g''(t)} {g
'(t)}|+(\frac {g''(t)} {g '(t)})^{2}+ |\frac {g'''(t)} {g '(t)}|)$
и
$$X_{r,g,\overline{x}}=\{ (q_1,...,q_n): \ q_1,...,q_n \in E_\delta $$
$$|\sum \limits_{k=1}\limits^{n}[(x_k - x_{k-1})(\frac
{g''(x_{k-1})} {g '(x_{k-1})} q_k'(0)- \frac {g ''(x_k)} {g
'(x_k)} q_k'(1))+$$
$$+(x_k - x_{k-1})^2\int_0^1 S_g (x_{k-1}+(x_k
- x_{k-1})q_k(t)) (q_k'(t))^2 dt ]|\leq 4 c_4 C_g r \}.
$$

{\bf Lemma  8.\/} {\it If $ \epsilon \in (0,1)$, then the
following inequality  is fulfilled \\ $\nu _n(E_\delta^n
\smallsetminus X_{\sqrt[3]{\epsilon},g,\overline{x}})\leq 2
\sqrt[3]{\epsilon}$ for any $g \in  \Diff^3_+([0,1])$,  for any
positive integer  $n$ and
 $\overline{x}=(x_1,...,x_{n-1}) \in D_n$, satisfying
 the inequality \\
 $\max \limits_{1\leq k \leq n}(x_k -
x_{k-1})<\epsilon $.\\
\/}

Proof. Let
$$ f_1 (q_1,...,q_n)=\sum \limits_{k=1}\limits^{n}(x_k - x_{k-1})(\frac
{g''(x_{k-1})} {g '(x_{k-1})} q_k'(0)- \frac {g ''(x_k)} {g
'(x_k)} q_k'(1)).$$

Then
$$
I_1= \int\limits_{E_\delta}...\int\limits_{E_\delta}f_1
(q_1,...,q_n)\,\nu (d q_1)... \nu (d q_n)=M_1 \sum
\limits_{k=1}\limits^{n}(x_k - x_{k-1})(\frac {g''(x_{k-1})} {g
'(x_{k-1})} - \frac {g ''(x_k)} {g '(x_k)} ).$$

As $ |\frac {g''(x_{k-1})} {g '(x_{k-1})} - \frac {g ''(x_k)} {g
'(x_k)}| \leq C_g (x_k - x_{k-1}),$ we have
$$
|I_1| \leq M_1 C_g  \sum \limits_{k=1}\limits^{n}(x_k - x_{k-1})^2
\leq M_1 C_g \epsilon \sum \limits_{k=1}\limits^{n}(x_k -
x_{k-1})=M_1 C_g \epsilon.$$

If $k\neq l$ then
$$\int\limits_{E_\delta}\int\limits_{E_\delta}
(\frac {g''(x_{k-1})} {g '(x_{k-1})}(q_k'(0)-M_1) - \frac {g
''(x_k)} {g '(x_k)} (q_k'(1)-M_1))$$
$$(\frac {g''(x_{l-1})} {g '(x_{l-1})}(q_l'(0)-M_1) - \frac {g
''(x_l)} {g '(x_l)}(q_l'(1)-M_1) )\nu (d q_k) \nu (d q_l)=0,$$
therefore
$$
I_2= \int\limits_{E_\delta}...\int\limits_{E_\delta}(f_1
(q_1,...,q_n)-I_1)^2\,\nu (d q_1)... \nu (d q_n)= $$
$$\sum \limits_{k=1}\limits^{n}(x_k - x_{k-1})^2\int\limits_{E_\delta}[\frac {g''(x_{k-1})} {g
'(x_{k-1})}(q_k'(0)-M_1) - \frac {g ''(x_k)} {g '(x_k)}
(q_k'(1)-M_1)]^2 \nu (d q_k)\leq$$
$$\leq 2 \sum \limits_{k=1}\limits^{n}(x_k - x_{k-1})^2[(\frac {g''(x_{k-1})} {g
'(x_{k-1})})^2 \int\limits_{E_\delta}(q_k'(0)-M_1)^2\nu (d q_k)
+$$
$$+(\frac {g ''(x_k)} {g '(x_k)})^2 \int\limits_{E_\delta}(q_k'(1)-M_1)^2
\nu (d q_k)]=$$
$$= 2 \sum \limits_{k=1}\limits^{n}(x_k - x_{k-1})^2[(\frac {g''(x_{k-1})} {g
'(x_{k-1})})^2 +(\frac {g ''(x_k)} {g '(x_k)})^2]
[\int\limits_{E_\delta}(q_k'(0))^2\nu (d q_k)-(M_1)^2] \leq$$
$$\leq 4 M_2 C_g \sum \limits_{k=1}\limits^{n}(x_k - x_{k-1})^2
 \leq 4 M_2 C_g \epsilon \sum \limits_{k=1}\limits^{n}(x_k
- x_{k-1})=4 M_2 C_g \epsilon.$$

Hence,
$$\nu _n ( \{ (q_1,...,q_n): | f_1
(q_1,...,q_n)-I_1| \geq 2 c_4 C_g \sqrt[3]{\epsilon} \} ) \leq
\frac{I_2}{(2 c_4 C_g \sqrt[3]{\epsilon})^2} \leq \frac{4 M_2 C_g
\epsilon}{(2 c_4 C_g \sqrt[3]{\epsilon})^2} \leq
\sqrt[3]{\epsilon}.
$$
Thus
$$\nu _n ( \{ (q_1,...,q_n): | f_1
(q_1,...,q_n)| \geq 3 c_4 C_g\sqrt[3]{\epsilon} \} ) \leq
\sqrt[3]{\epsilon}.
$$

Let
$$ f_2 (q_1,...,q_n)=\sum \limits_{k=1}\limits^{n}((x_k - x_{k-1})^2\int_0^1 S_g (x_{k-1}+(x_k
- x_{k-1})q_k(t)) (q_k'(t))^2 dt.$$

Then
$$
I_3= \int\limits_{E_\delta}...\int\limits_{E_\delta}|f_2
(q_1,...,q_n)|\,\nu (d q_1)... \nu (d q_n)\leq $$
$$\leq 2 C_g \sum
\limits_{k=1}\limits^{n}(x_k - x_{k-1})^2 \int\limits_{E_\delta}(
\int_0^1 (q_k'(t))^2 dt)\nu (d q_k)\leq$$
$$\leq 2 c_4 C_g \epsilon \sum
\limits_{k=1}\limits^{n}(x_k - x_{k-1})=2 c_4 C_g \epsilon .$$
Thus
$$\nu _n ( \{ (q_1,...,q_n): | f_2
(q_1,...,q_n)| \geq 2 c_4 C_g\sqrt[3]{\epsilon} \} ) \leq
\frac{I_3}{2 c_4 C_g \sqrt[3]{\epsilon}} \leq \frac{2 c_4 C_g
\epsilon}{2 c_4 C_g \sqrt[3]{\epsilon}}= (\sqrt[3]{\epsilon})^2
\leq \sqrt[3]{\epsilon}.
$$

Hence,
$$\nu _n(E_\delta ^n \smallsetminus
X_{\sqrt[3]{\epsilon},g,\overline{x}})= $$
$$=\nu _n ( \{ (q_1,...,q_n): |f_1
(q_1,...,q_n)+ f_2 (q_1,...,q_n)| \geq 4 c_4 C_g\sqrt[3]{\epsilon}
\} )\leq $$
$$\leq \nu _n ( \{ (q_1,...,q_n): |f_1
(q_1,...,q_n)| \geq 2 c_4 C_g\sqrt[3]{\epsilon} \} )+$$
$$+\nu _n ( \{ (q_1,...,q_n): |f_2 (q_1,...,q_n)| \geq 2 c_4 C_g\sqrt[3]{\epsilon}
\} )\leq 2 \sqrt[3]{\epsilon},$$   which implies the assertion of
Lemma 8.

{\bf Lemma 9.\/} {\it For any $g \in  \Diff^3_0([0,1])$,
 $ \epsilon>0$,  there are $r>1$ and
$\delta_1 \in (0,1)$ such that the inequality is valid
$$|\prod \limits_{k=1}
\limits^{n}\frac{v(\frac{g(x_k)-g(x_{k-2})}{2\sqrt{(g(x_k)-g(x_{k-1}))(g(x_{k-1})-g(x_{k-2}))}})}
{v(\frac{x_k-x_{k-2}}{2\sqrt{(x_k-x_{k-1})(x_{k-1}-x_{k-2})}})}-1|\leq
\epsilon,$$ $g(x_{-1})=g(x_{n-1})-1$ for any natural $n$ and any
$\overline{x}=(x_1,...,x_{n-1}) \in D_n$ satisfying the inequality
$\max \limits_{1\leq k \leq n}(x_k - x_{k-1})<\delta_1$, $\min
\limits_{1\leq k \leq
n}(\frac{x_k-x_{k-2}}{2\sqrt{(x_k-x_{k-1})(x_{k-1}-x_{k-2})}})>r$,
where  $x_0=0$,  $x_n=1$, $x_{-1}=x_{n-1}-1$,.\/}

Proof. Let $\epsilon \in (0,1)$. Let $C=\max \limits_{t_1, t_2 \in
[0,1]}|\frac{g''(t_{1})}{g'(t_{2})}|$, $\delta_1
=\frac{1}{400(C+1)}$ и $r=e^{\frac{8000 (C+1)}{\epsilon}} $.

For any $k \ (1\leq k \leq n)$, there are $x_{k-1}^*, x_{k-1}^{**}
\in (0,1)$ such that
$$g(x_k)-g(x_{k-1})=g'(x_{k-1})(x_k-x_{k-1})+\frac{1}{2}
g''(x_{k-1}^*)(x_k-x_{k-1})^2=$$
$$=g'(x_{k-1})(x_k-x_{k-1})(1+\frac{g''(x_{k-1}^*)}{2g'(x_{k-1})}
(x_k-x_{k-1})),$$
$$g(x_{k-1})-g(x_{k-2})=g'(x_{k-1})(x_{k-1}-x_{k-2})+\frac{1}{2}
g''(x_{k-1}^{**})(x_{k-1}-x_{k-2})^2$$
$$=g'(x_{k-1})(x_{k-1}-x_{k-2})(1+\frac{g''(x_{k-1}^{**})}{2g'(x_{k-1})}
(x_{k-1}-x_{k-2})),$$ for $k\geq 2$ and for $k=1$
$$g(x_1)-g(x_0)=g'(x_0)(x_1-x_0)+\frac{1}{2}
g''(x_0^*)(x_1-x_0)^2=$$
$$=g'(x_0)(x_1-x_0)(1+\frac{g''(x_0^*)}{2g'(x_0)}
(x_1-x_0)),$$
$$g(x_0)-g(x_{-1})=g(1)-g(x_{n-1})=g'(1)(1-x_{n-1})+\frac{1}{2}
g''(x_0^{**})(1-x_{n-1})^2=$$
$$=g'(x_0)(x_0-x_{-1})(1+\frac{g''(x_0^{**})}{2g'(x_0)}
(x_0-x_{-1})).$$

Hence
$$\frac{g(x_k)-g(x_{k-2})}{2\sqrt{(g(x_k)-g(x_{k-1}))(g(x_{k-1})-g(x_{k-2}))}}=$$
$$=\frac{x_k-x_{k-2}}{2\sqrt{(x_k-x_{k-1})(x_{k-1}-x_{k-2})}}
(1+\lambda _k (x_k-x_{k-2})),
$$
where
$$\lambda _k =\frac
{1}{x_k-x_{k-2}} (\frac{1+
\frac{g''(x_{k-1}^*)(x_k-x_{k-1})^2+g''(x_{k-1}^{**})(x_{k-1}-x_{k-2})^2}
{2g'(x_{k-1})(x_{k}-x_{k-2})}}
{\sqrt{(1+\frac{g''(x_{k-1}^*)}{2g'(x_{k-1})}
(x_k-x_{k-1}))(1+\frac{g''(x_{k-1}^{**})}{2g'(x_{k-1})}
(x_{k-1}-x_{k-2}))}}-1).
$$

As $C(x_k-x_{k-2})<C\delta_1 =\frac{1}{400}$, then $|\lambda _k|
<4C$.

Let $t=\frac{x_k-x_{k-2}}{2\sqrt{(x_k-x_{k-1})(x_{k-1}-x_{k-2})}}$
and $\alpha = \lambda _k (x_k-x_{k-2}).$ Then $$t(1+\alpha )
=\frac{g(x_k)-g(x_{k-2})}{2\sqrt{(g(x_k)-g(x_{k-1}))(g(x_{k-1})-g(x_{k-2}))}}$$
and $|\alpha |<\frac{1}{80}$,  $t>4 $.

Let us obtain $\tau = \ln(t+\sqrt{t^2-1})> \ln t>1 $       and
$$\beta = \ln [1+\alpha (\frac{t}{t+\sqrt{t^2-1}})(1+\frac{(2+\alpha ) t}
{\sqrt{t^2-1}+\sqrt{t^2(1+\alpha )^2-1}})].$$ Then
$\frac{1}{4}|\alpha|\leq |\beta |\leq 4|\alpha|$, $ |\beta |\leq
\frac{1}{20}$ and
$$\ln(t(1+\alpha )+\sqrt{t^2(1+\alpha )^2-1})=\ln(t+\sqrt{t^2-1})
+
$$
$$ \ln [1+\alpha (\frac{t}{t+\sqrt{t^2-1}})(1+\frac{(2+\alpha ) t}
{\sqrt{t^2-1}+\sqrt{t^2(1+\alpha )^2-1}})]=\tau+\beta.
$$
There exists $ \theta \in (0,1)$ such that $v(t(1+\alpha
))=v_1(\tau+\beta)=v_1(\tau)+v_1'(\tau+\theta \beta)\beta.
$

It follows from the lemma 2 that $ |v_1'(\tau+\theta \beta)|\leq
\frac{4}{\tau+\theta \beta}v_1(\tau+\theta \beta). $ As $|\theta
\beta|\leq \frac{1}{20}$, we have $ |v_1'(\tau+\theta \beta)|\leq
\frac{10}{\tau}v_1(\tau). $

Thus it exists $\omega (\alpha ,t) $ such that $|\omega (\alpha
,t)|<10$ и $\frac{v_1(\tau+\beta)}{v_1(\tau)}=1+ \omega (\alpha
,t)\frac{\beta}{\tau}.$ Taking account  $\frac{1}{4}|\lambda _k
(x_k-x_{k-2})|=\frac{1}{4}|\alpha|\leq |\beta |\leq
4|\alpha|=4|\lambda _k (x_k-x_{k-2})|$, we receive
$$\frac{v(\frac{g(x_k)-g(x_{k-2})}{2\sqrt{(g(x_k)-g(x_{k-1}))(g(x_{k-1})-g(x_{k-2}))}})}
{v(\frac{x_k-x_{k-2}}{2\sqrt{(x_k-x_{k-1})(x_{k-1}-x_{k-2})}})}=$$
$$=\frac{v_1(\tau+\beta)}{v_1(\tau)}=1+ \omega _k
\frac{x_k-x_{k-2}}{\ln
\frac{x_k-x_{k-2}}{2\sqrt{(x_k-x_{k-1})(x_{k-1}-x_{k-2})}}},
$$
where $\omega _k= \omega (\alpha ,t)\frac{\beta \ln t}{\tau
(x_k-x_{k-2})}$ и $|\omega _k|\leq 200 C $.

Hence
$$\sigma=\ln(\prod \limits_{k=1}
\limits^{n}\frac{v(\frac{g(x_k)-g(x_{k-2})}{2\sqrt{(g(x_k)-g(x_{k-1}))(g(x_{k-1})-g(x_{k-2}))}})}
{v(\frac{x_k-x_{k-2}}{2\sqrt{(x_k-x_{k-1})(x_{k-1}-x_{k-2})}}))=}
$$
$$=\sum \limits_{k=1}
\limits^{n}\ln(1+ \omega _k \frac{x_k-x_{k-2}}{\ln
\frac{x_k-x_{k-2}}{2\sqrt{(x_k-x_{k-1})(x_{k-1}-x_{k-2})}}})
$$
and
$$|\sigma | \leq 2\sum \limits_{k=1}
\limits^{n}|\omega _k |\frac{x_k-x_{k-2}}{\ln
\frac{x_k-x_{k-2}}{2\sqrt{(x_k-x_{k-1})(x_{k-1}-x_{k-2})}}}\leq
\frac{400C}{\ln r}\sum \limits_{k=1} \limits^{n}{(x_k-x_{k-2})}=
\frac{800C}{\ln r}\leq \frac{\epsilon}{10},
$$
therefore
$$|\prod \limits_{k=1}
\limits^{n}\frac{v(\frac{g(x_k)-g(x_{k-2})}{2\sqrt{(g(x_k)-g(x_{k-1}))(g(x_{k-1})-g(x_{k-2}))}})}
{v(\frac{x_k-x_{k-2}}{2\sqrt{(x_k-x_{k-1})(x_{k-1}-x_{k-2})}})}-1|=
$$
$$=|e^ \sigma   -1| \leq  e^{\frac{\epsilon}{10}}-e^{-\frac{\epsilon}{10}}
\leq  \frac{\epsilon}{5}+\frac{\epsilon}{5} < \epsilon,
$$
which implies the assertion of Lemma 9.

Introduce the mapping $Q_n:D_n \times E_\delta^n \to
E_\delta=\Diff^{1,\delta}_+([0,1])$ by setting \\
$f_n\circ
(l_n)^{-1}
 =Q_n(x_1,...,x_{n-1},\varphi_1,...,\varphi_{n})$, where
$$f_n(t)=x_{k-1}+(x_k-x_{k-1})\varphi_{k}(n(t-\frac{k-1}{n})),$$
$$l_n (t)=\frac{1}
{x_1-x_0+\sum \limits_{m=2}
\limits^{n}(x_m-x_{m-1})\frac{\varphi_2'(0)\varphi_3'(0)...\varphi_{m}'(0)}
{\varphi_1'(1)\varphi_2'(1)...\varphi_{m-1}'(1)}}\cdot$$
$$\cdot(x_1-x_0+\sum
\limits_{m=2}\limits^{k-1}(x_m-x_{m-1})\frac{\varphi_2'(0)\varphi_3'(0)...\varphi_{m}'(0)}
{\varphi_1'(1)\varphi_2'(1)...\varphi_{m-1}'(1)}+$$
$$+(x_k-x_{k-1})\frac{\varphi_2'(0)\varphi_3'(0)...\varphi_{k}'(0)}
{\varphi_1'(1)\varphi_2'(1)...\varphi_{k-1}'(1)}n(t-\frac{k-1}{n}))
$$
for $t\in [\frac{k-1}{n},\frac{k}{n}]$, $(x_1,...,x_{n-1})\in
D_{n}$, $(\varphi_1,...,\varphi_{n})\in E_\delta^n$.

The function $f=f_n\circ (l_n)^{-1}$ belongs to
$\Diff^{1,\delta}_+([0,1])$, because the left derivation
$$f'( (l_n)^{-1}(\frac{k-1}{n}-0))=n(x_{k-1}-x_{k-2})\varphi_{k-1}(1)\cdot$$
$$\cdot\frac {x_1-x_0+\sum \limits_{m=2}
\limits^{n}(x_m-x_{m-1})\frac{\varphi_2'(0)\varphi_3'(0)...\varphi_{m}'(0)}
{\varphi_1'(1)\varphi_2'(1)...\varphi_{m-1}'(1)}}
{(x_{k-1}-x_{k-2})\frac{\varphi_2'(0)\varphi_3'(0)...\varphi_{k-1}'(0)}
{\varphi_1'(1)\varphi_2'(1)...\varphi_{k-2}'(1)}n }=$$
$$= (x_1-x_0+\sum \limits_{m=2}
\limits^{n}(x_m-x_{m-1})\frac{\varphi_2'(0)\varphi_3'(0)...\varphi_{m}'(0)}
{\varphi_1'(1)\varphi_2'(1)...\varphi_{m-1}'(1)}) \frac
{\varphi_1'(1)\varphi_2'(1)...\varphi_{k-1}'(1)}
{\varphi_2'(0)\varphi_3'(0)...\varphi_{k-1}'(0)}.$$ is equal to
the right derivation
$$f'( (l_n)^{-1}(\frac{k-1}{n}+0))=n(x_k-x_{k-1})\varphi_{k}(0)\cdot$$
$$\cdot\frac {x_1-x_0+\sum \limits_{m=2}
\limits^{n}(x_m-x_{m-1})\frac{\varphi_2'(0)\varphi_3'(0)...\varphi_{m}'(0)}
{\varphi_1'(1)\varphi_2'(1)...\varphi_{m-1}'(1)}}
{(x_k-x_{k-1})\frac{\varphi_2'(0)\varphi_3'(0)...\varphi_{k}'(0)}
{\varphi_1'(1)\varphi_2'(1)...\varphi_{k-1}'(1)}n }=$$
$$= (x_1-x_0+\sum \limits_{m=2}
\limits^{n}(x_m-x_{m-1})\frac{\varphi_2'(0)\varphi_3'(0)...\varphi_{m}'(0)}
{\varphi_1'(1)\varphi_2'(1)...\varphi_{m-1}'(1)}) \frac
{\varphi_1'(1)\varphi_2'(1)...\varphi_{k-1}'(1)}
{\varphi_2'(0)\varphi_3'(0)...\varphi_{k-1}'(0)}.$$

Write $$L_{\delta , n}(F)=\int\limits_{0}\limits^{1}
\int\limits_{x_1}\limits^{1}...\int\limits_{x_{n-2}}\limits^{1}
\int\limits_{E_\delta}...\int\limits_{E_\delta}
F(Q_n(x_1,x_2,...,x_{n-1},\varphi_1,...,\varphi_{n}))\cdot$$
$$\cdot u(x_1,x_2,...,x_{n-1})dx_1dx_2...dx_{n-1}\nu (d\varphi_1)...\nu (d\varphi_{n})$$
for any function $F \in
C_b(E_\delta)=C_b(\Diff^{1,\delta}_+([0,1]))$.

{\bf Theorem 3.\/} {\it  a function $F $ belongs to the space
$C_b(\Diff^{1,\delta}_+([0,1]))$, let \\ a diffeomorphism  $g$
belongs to the group $\Diff^3_0([0,1])$.\\ Then $\lim \limits_{n
\to \infty}|L_{\delta , n}(F_g )-L_{\delta , n}(F)|=0$. \/}

Proof. Let $F \in C_b(\Diff^{1,\delta}_+([0,1]))$, $C= \sup
\limits_{g \in E_\delta}|F(f)|.$

 Let $\epsilon
\in (0,1)$ and $g \in \Diff^3_0([0,1])$. Let us take positive
$\epsilon _1$ satisfying the inequalities $\epsilon
_1<\frac{1}{8}\epsilon$, \\
$e^{4c_5 C_g \sqrt[3]{\epsilon _1} }-e^{-4c_4 C_g
\sqrt[3]{\epsilon _1} }<\epsilon$.

It follows from the lemma 9 that there are $r>1$ and $\delta _1
\in (0,1)$ such that
$$|\prod \limits_{k=1}
\limits^{n}\frac{v(\frac{g(x_k)-g(x_{k-2})}{2\sqrt{(g(x_k)-g(x_{k-1}))(g(x_{k-1})-g(x_{k-2}))}})}
{v(\frac{x_k-x_{k-2}}{2\sqrt{(x_k-x_{k-1})(x_{k-1}-x_{k-2})}})}-1|\leq
\epsilon$$ for any positive integer  $n$ and for
$\overline{x}=(x_1,...,x_{n-1}) \in D_n$ satisfying the
inequalities $\max \limits_{1\leq k \leq n}(x_k -
x_{k-1})<\delta_1$, $\min \limits_{1\leq k \leq
n}(\frac{x_k-x_{k-2}}{2\sqrt{(x_k-x_{k-1})(x_{k-1}-x_{k-2})}})>r$.

 It follows from the lemma 8 that the inequality is valid $\nu _n(E_\delta ^n
\smallsetminus X_{\sqrt[3]{\epsilon_1},g,\overline{x}})\leq 2
\sqrt[3]{\epsilon_1}\leq \epsilon$ for any positive integer $n$
and for any $\overline{x}=(x_1,...,x_{n-1}) \in D_n$ satisfying
the inequalities   $\max \limits_{1\leq k \leq n}(x_k -
x_{k-1})<\epsilon _1$.

It follows from the lemma 5, 6 that there exists a positive
integer $N$ such that
$$\int\limits_{0}\limits^{1} \int\limits_{x_1}\limits^{1}...\int\limits_{x_{n-2}}\limits^{1}
(1-\vartheta (\frac{1}{\min(\epsilon_1 ,\delta )}\max
\limits_{1\leq k \leq n} (x_k-x_{k-1})))
u_{n}(x_1,x_2,...,x_{n-1})dx_1dx_2...dx_{n-1}<\epsilon ,$$
$$\int\limits_{0}\limits^{1} \int\limits_{x_1}\limits^{1}...\int\limits_{x_{n-2}}\limits^{1}
\vartheta (\frac{1}{r}\min \limits_{1\leq k \leq n}
\frac{x_{k}-x_{k-2}}{2\sqrt{(x_{k}-x_{k-1})(x_{k-1}-x_{k-2})}})$$
$$u_{n}(x_1,x_2,...,x_{n-1})dx_1dx_2...dx_{n-1}<\epsilon $$
for any integer $n>N$.

 Let
$$\phi (x_1,x_2,...,x_{n-1})=\vartheta
(\frac{1}{\min(\epsilon_1 ,\delta )}\max \limits_{1\leq k \leq n}
(x_k-x_{k-1})))$$
$$(1-\vartheta (\frac{1}{r}\min \limits_{1\leq k
\leq n}
\frac{x_{k}-x_{k-2}}{2\sqrt{(x_{k}-x_{k-1})(x_{k-1}-x_{k-2})}})),$$
then
$$\int\limits_{0}\limits^{1} \int\limits_{x_1}\limits^{1}...\int\limits_{x_{n-2}}\limits^{1}
\phi (x_1,x_2,...,x_{n-1})
u_{n}(x_1,x_2,...,x_{n-1})dx_1dx_2...dx_{n-1}\geq 1-$$
$$-\int\limits_{0}\limits^{1} \int\limits_{x_1}\limits^{1}...\int\limits_{x_{n-2}}\limits^{1}
(1-\vartheta (\frac{1}{\min(\epsilon _1 ,\delta )}\max
\limits_{1\leq k \leq n} (x_k-x_{k-1})))
u_{n}(x_1,x_2,...,x_{n-1})dx_1dx_2...dx_{n-1}-$$
$$-\int\limits_{0}\limits^{1} \int\limits_{x_1}\limits^{1}...\int\limits_{x_{n-2}}\limits^{1}
\vartheta (\frac{1}{r}\min \limits_{1\leq k \leq n}
\frac{x_{k}-x_{k-2}}{2\sqrt{(x_{k}-x_{k-1})(x_{k-1}-x_{k-2})}})$$
$$u_{n},x_2,...,x_{n-1})dx_1dx_2...dx_{n-1}>1-2\epsilon .$$

Let $y_k=g(x_k)$, $\overline{y}=(y_1,...,y_{n-1})\in D_n$, \\
$g^{-1}(\overline{y})=(x_1,...,x_{n-1})\in D_n$, $\phi _g
(y_1,...,y_{n-1})=\phi (x_1,...,x_{n-1}) $ and
$$gX_{\sqrt[3]{\epsilon_1},g,g^{-1}(\overline{y})}=\{(\varphi_1,...,\varphi_{n})
: (q_1,...,q_{n-1})\in
 X_{\sqrt[3]{\epsilon_1},g,\overline{x}}, \ \ $$
$$ Q_n(y_1,...,y_{n-1},\varphi_1,...,\varphi_{n})=g\circ
(Q_n(x_1,x_2,...,x_{n-1},q_1,...,q_{n-1}))\}
 .$$

It is easy to see that
$\varphi_{k}(t)=\frac{g(x_{k-1}+(x_{k}-x_{k-1})q_k
(nt-k+1))-g(x_{k-1})}{g(x_{k})-g(x_{k-1})}, $ because
$$\frac{(y_k-y_{k-1})\varphi_2'(0)\varphi_3'(0)...\varphi_{k}'(0)}
{(y_1-y_{0})\varphi_1'(1)\varphi_2'(1)...\varphi_{k-1}'(1)}=
\frac{(x_k-x_{k-1})q_2'(0)q_3'(0)...q_{k}'(0)}
{(x_1-x_{0})q_1'(1)q_2'(1)...q_{k-1}'(1)}$$

We have
$$\int\limits_{0}\limits^{1} \int\limits_{y_1}\limits^{1}...\int\limits_{y_{n-2}}\limits^{1}
\phi _g  (\overline{y}) u_{n}(\overline{y}) \nu
_n(gX_{\sqrt[3]{\epsilon_1},g,g^{-1}(\overline{y})})
dy_1dy_2...dy_{n-1}=
$$
$$=\int\limits_{0}\limits^{1} \int\limits_{x_1}\limits^{1}...\int\limits_{x_{n-2}}\limits^{1}
\phi (\overline{x}) u_{n}(\overline{x})
(\int\limits_{X_{\sqrt[3]{\epsilon_1},g,\overline{x}}} exp(\sum
\limits_{k=1}\limits^{n}[(x_k - x_{k-1})(\frac {g''(x_{k-1})} {g
'(x_{k-1})} q_k'(0)- \frac {g ''(x_k)} {g '(x_k)} q_k'(1))+$$
$$+(x_k - x_{k-1})^2\int_0^1 S_g (x_{k-1}+(x_k
- x_{k-1})q_k(t)) (q_k'(t))^2 dt ])\nu (d q_1)... \nu (d q_n))
$$
$$\prod \limits_{k=1}
\limits^{n}\frac{v(\frac{g(x_k)-g(x_{k-2})}{2\sqrt{(g(x_k)-g(x_{k-1}))(g(x_{k-1})-g(x_{k-2}))}})}
{v(\frac{x_k-x_{k-2}}{2\sqrt{(x_k-x_{k-1})(x_{k-1}-x_{k-2})}})}
dx_1dx_2...dx_{n-1}\geq $$
$$\geq(1-\epsilon)^2\int\limits_{0}\limits^{1}
\int\limits_{x_1}\limits^{1}...\int\limits_{x_{n-2}}\limits^{1}
\phi (\overline{x}) u_{n}(\overline{x})\nu
_n(X_{\sqrt[3]{\epsilon_1},g,\overline{x}})dx_1dx_2...dx_{n-1}\geq$$
$$\geq(1-\epsilon)^3\int\limits_{0}\limits^{1}
\int\limits_{x_1}\limits^{1}...\int\limits_{x_{n-2}}\limits^{1}
\phi (\overline{x}) u_{n}(\overline{x})dx_1dx_2...dx_{n-1}\geq
(1-\epsilon)^3(1-2\epsilon).$$

Hence,
$$|L_{\delta ,n} (F_g )-\int\limits_{0}\limits^{1} \int\limits_{y_1}\limits^{1}...\int\limits_{y_{n-2}}\limits^{1}
\phi _g  (\overline{y}) u_{n}(\overline{y})(
\int\limits_{gX_{\sqrt[3]{\epsilon_1},g,g^{-1}(\overline{y})}}F_g
(Q_n(\overline{y},\varphi_1,...,\varphi_{n}))$$
$$\nu (d\varphi_1)...\nu (d\varphi_{n}))
dy_1dy_2...dy_{n-1}|\leq C(1-(1-\epsilon)^3(1-2\epsilon))
$$
and
$$|L_{\delta ,n} (F)-\int\limits_{0}\limits^{1} \int\limits_{x_1}\limits^{1}...\int\limits_{x_{n-2}}\limits^{1}
\phi _g  (\overline{x}) u_{n}(\overline{x})(
\int\limits_{X_{\sqrt[3]{\epsilon_1},g,\overline{x}}}
F(Q_n(\overline{x},q_1,...,q_{n}))$$
$$\nu (dq_1)...\nu (dq_{n}))
dx_1dx_2...dx_{n-1}|\leq C(1-(1-\epsilon)(1-2\epsilon)).
$$

We have
$$|\int\limits_{0}\limits^{1} \int\limits_{y_1}\limits^{1}...\int\limits_{y_{n-2}}\limits^{1}
\phi _g  (\overline{y}) u_{n}(\overline{y})(
\int\limits_{gX_{\sqrt[3]{\epsilon_1},g,g^{-1}(\overline{y})}}F_g
(Q_n(\overline{y},\varphi_1,...,\varphi_{n}))$$
$$\nu (d\varphi_1)...\nu (d\varphi_{n}))
dy_1dy_2...dy_{n-1}-
$$
$$-\int\limits_{0}\limits^{1} \int\limits_{x_1}\limits^{1}...\int\limits_{x_{n-2}}\limits^{1}
\phi _g  (\overline{x}) u_{n}(\overline{x})(
\int\limits_{X_{\sqrt[3]{\epsilon_1},g,\overline{x}}}
F(Q_n(\overline{x},q_1,...,q_{n}))$$
$$\nu (dq_1)...\nu (dq_{n}))
dx_1dx_2...dx_{n-1}|\leq
$$
$$\leq\int\limits_{0}\limits^{1} \int\limits_{x_1}\limits^{1}...\int\limits_{x_{n-2}}\limits^{1}
\phi (\overline{x}) u_{n}(\overline{x})
(\int\limits_{X_{\sqrt[3]{\epsilon_1},g,\overline{x}}} |exp(\sum
\limits_{k=1}\limits^{n}[(x_k - x_{k-1})(\frac {g''(x_{k-1})} {g
'(x_{k-1})} q_k'(0)- \frac {g ''(x_k)} {g '(x_k)} q_k'(1))+$$
$$+(x_k - x_{k-1})^2\int_0^1 S_g (x_{k-1}+(x_k
- x_{k-1})q_k(t)) (q_k'(t))^2 dt ]
$$
$$\prod \limits_{k=1}
\limits^{n}\frac{v(\frac{g(x_k)-g(x_{k-2})}{2\sqrt{(g(x_k)-g(x_{k-1}))(g(x_{k-1})-g(x_{k-2}))}})}
{v(\frac{x_k-x_{k-2}}{2\sqrt{(x_k-x_{k-1})(x_{k-1}-x_{k-2})}})}-1|$$
$$|F(Q_n(\overline{x},q_1,...,q_{n}))|\nu (d q_1)... \nu (d
q_n)) u_{n}(\overline{x})dx_1dx_2...dx_{n-1}\leq
$$
$$\leq C\epsilon(2+\epsilon)\int\limits_{0}\limits^{1}
\int\limits_{x_1}\limits^{1}...\int\limits_{x_{n-2}}\limits^{1}
\phi (\overline{x}) u_{n}(\overline{x})\nu
_n(X_{\sqrt[3]{\epsilon_1},g,\overline{x}})dx_1dx_2...dx_{n-1}\leq
C\epsilon(2+\epsilon),$$ which implies the assertion of Theorem 3.

Define a ultrafilter  $\Im$ on the set positive integers such that
 $\Im$ contains the sets $ \{n,n+1,...\}$ for any positive integer $n$.
We set $L_\delta(F)=\lim \limits_{\Im}L_{\delta , n} ( F)$ for any
function $F \in C_b(E_\delta)$.

Note that the limit always exists because $|L_{\delta , n} (
F)|\leq \sup \limits_{f \in E_\delta}| F(f)|.$

It is easy to see that $L(e_{1,\delta})=1$, $|L_{\delta } (
F)|\leq \sup \limits_{f \in E_\delta}| F(f)|,$ and $L(F)\geq 0$
for any nonnegative function $F \in B(G)$. In turn, Theorem 1
follows from Theorem 3.

{\bf Acknowledgements.} The author acknowledge I.K. Babenko, R.
Grigorchuk, P. de la Harpe, V.Sergiescu, O.G. Smolyanov, A.I.
Stern for the discussions and interest about this work.

 {\bf References}

1. Cannon J.W., Floyd W.J., Parry W.R. "Introductory notes on
Richard Thompson's groups", Enseign Math., vol 42, issue 2 (1996),
pages 215--256.

2. Ghys \`E., Sergiescu V. "Sur un groupe remarquable de
diffeomorphismes du cercle", Comment.Math.Helvetici 62 (1987)
185--239.

3.Shavgulidze E.T. "Some Properties of Quasi-Invariant Measures on
Groups of Diffeomorphisms of the Circle",  Russ. J. Math. Physics,
vol 7, issue 4 (2000 ), pages 464--472.

4.Shavgulidze E.T. "Amenability of discrete subgroups of the group
of diffeomorphisms of the circle", Russ. J. Math. Physics, vol 16,
issue 1 (2009 ), pages 138--140

5. Hui-Hsiuhg Kuo "Gaussian measures in Banach space", Lecture
Notes in Mathematics 463, Springer-Verlag 1975.

6. G.H.Watson "A Treatise on the Theory of Bessel Functions",
Cambridge 1944.

\end{document}